\documentclass[twoside,british,english]{elsarticle}
\usepackage[T1]{fontenc}
\usepackage[latin9]{inputenc}
\usepackage{array}
\usepackage{booktabs}
\usepackage{url}
\usepackage{multirow}
\usepackage{amsmath}
\usepackage{amssymb}
\usepackage{graphicx}
\usepackage{esint}

\usepackage[labelfont=bf,justification=raggedright,singlelinecheck=false]{caption}
\makeatletter

\newcommand{\CC}{C\nolinebreak\hspace{-.05em}\raisebox{.4ex}{\tiny\bf +}\nolinebreak\hspace{-.10em}\raisebox{.4ex}{\tiny\bf +}}

\providecommand{\tabularnewline}{\\}

\journal{Applied Mathematics and Computation}

\usepackage{booktabs}

\makeatother

\usepackage{babel}
\begin{document}

\title{PoisFFT - A Free Parallel Fast Poisson Solver\tnoteref{t1}}

\tnotetext[PoisFFT]{The source code for this software can be found at \url{https://github.com/LadaF/PoisFFT}.}

\author[KMOP]{V.~Fuka\corref{cor1}}

\ead{vladimir.fuka@mff.cuni.cz}

\cortext[cor1]{Corresponding author}

\address[KMOP]{Department of Meteorology and Environment Protection, Faculty of
Mathematics and Physics, V Holešovi\v{c}kách 2, Prague 8, Czech Republic
18000}
\begin{abstract}
A fast Poisson solver software package PoisFFT is presented. It is
available as a free software licensed under the GNU GPL license version 3. The
package uses the fast Fourier transform to directly solve the Poisson
equation on a uniform orthogonal grid. It can solve the pseudo-spectral
approximation and the second order finite difference approximation
of the continuous solution. The paper reviews the mathematical methods
for the fast Poisson solver and discusses the software implementation
and parallelization. The use of PoisFFT in an incompressible flow
solver is also demonstrated.\end{abstract}
\begin{keyword}
Poisson equation \sep fast Poisson solver \sep parallel \sep MPI
\sep OpenMP
\end{keyword}
\maketitle

\section{Introduction}

The Poisson equation is an elliptic partial differential equation
with a wide range of applications in physics or engineering. In this
paper the form

\begin{align}
\nabla^{2}\varphi & =g,\label{eq:poisson-continuous}
\end{align}
in a rectangular domain $\Omega\subset\mathbb{R}^{d}$, for $d=1\ldots3$,
with boundary conditions specified on $\partial\Omega$ is considered.
The solution $\varphi$ and the specified source term $g$ are both
sufficiently smooth real functions on $\Omega$. The equation can
be used for computations of electrostatic and gravitational potentials,
to solve potential flow of ideal fluids or to perform pressure correction
when solving the incompressible Navier-Stokes equations.

There are many methods for the approximate numerical solution of the
Poisson equation. It can be discretized in many ways using the finite
difference, finite volume, finite element or (\foreignlanguage{british}{pseudo-})spectral
methods.

In many problems the goal is to compute the most accurate approximation
of the continuous problem, while in other applications it is necessary
to use the discretization consistent with other parts of a larger
solver. This is true for those incompressible flow solvers which enforce
the discrete continuity equation by solving the discrete
Poisson equation\citep{brown:proj}.

The fast direct solvers of the Poisson equation have a long history.
The method of cyclic reduction was used by \citet{Hockney:fast-direct},
\citet{BuzbeeGolubNielson:cyclic} and \citet{Buneman:cyclic}. The
Fourier analysis method was developed by \citet{CooleyLewisWelch},
\citet{WilhelmsonEricksen} and \citet{Swarztrauber:review, swarztrauber:1986symmetric}
for various boundary conditions. The methods of the cyclic reduction
and Fourier analysis can be combined in the FACR algorithm \citep{swarztrauber:1986symmetric}.
The discrete Fourier transforms used with different boundary conditions
were surveyed by \citet{SchumannSweet:FFTs}.

The requirement of a rectangular grid or other simple type of grid
can be too restrictive for many applications, but the fast solution
on a regular grid can be often used as a preconditioner for a solution
on an unstructured grid or for problems with non-constant coefficients
(e.g., $\nabla\cdot(a(x)\nabla\varphi)=g$)\citep{Gockenbach}.

Many implementations of fast Poisson solvers exist, but their software
implementation is often not publically available, or it is a part of a larger software
package and is able to solve only those cases, which are applicable
in that particular area \citep{laizet2011IncpowtoototacturprouptoO105comcor,Garcia-Risueno:ParallelPoisson}. One
of the more general ones is the package FISHPACK\citep{FISHPACK}
which solves the second order finite difference approximation of the
Poisson or Helmholtz equation on rectangular grids. Its origin predates
FORTRAN 77 and although it has been updated in 2004 to version 5.0
to be compatible with Fortran 90 compilers, the non-existent parallel
version, the remaining non-standard features and missing interfaces
to other programming languages make it obsolete.

This paper presents the free software implementation of a fast direct
solver for the Poisson equation in 1, 2 and 3 dimensions with several
types of boundary conditions and with both the pseudo-spectral and
the finite difference approximation. The solver is called PoisFFT.
Section \ref{sec:Solution-methods} contains the overview of the methods
used to solve the Poisson equations in PoisFFT. Section \ref{sec:Software-implementation}
briefly introduces the software implementation and its performance
is tested in Section \ref{sec:Performance-evaluation}. In Section
\ref{sec:Application-to-incompressible} we present the application
of the Poisson solver to the solution of the incompressible flow.

\section{Solution methods\label{sec:Solution-methods}}

\subsection{Pseudo-spectral method\label{sub:Pseudo-spectral-method}}

The first type of approximation PoisFFT solves is the pseudo-spectral
method approximation which is achieved by decomposing the continuous
problem using the Fourier series and solving it in the frequency domain.
Consider the 1D Poisson problem \eqref{eq:poisson-continuous} on
interval $[0,L]$. With periodic boundary conditions any sufficiently
smooth function $f$ on this interval can be decomposed as

\begin{align}
f(x) & =\sum_{k=-\infty}^{\infty}\widehat{f}_{k}e^{2\pi ikx/L},\label{eq:Four-series}
\end{align}
where the complex coefficients are computed as

\begin{align}
\hat{f}_{k} & =\frac{1}{L}\int_{0}^{L}f(x)e^{-2\pi ikx/L}dx.\label{eq:Four-series-coef}
\end{align}

After substituting from \eqref{eq:Four-series} for $\varphi$ and
$g$ we get

\begin{align}
-\left(\frac{2\pi k}{L}\right)^{2}\hat{\varphi}_{k} & =\hat{g}_{k,}\quad\forall k,\label{eq:spectral-1d-solution}
\end{align}
which can be used to compute the solution in the basis of the Fourier
modes. The modes $e^{2\pi ik_{x}x/L_{x}}$ are eigenvectors of the
linear operator $\frac{\partial^{2}}{\partial x^{2}}$ with the eigenvalues
$\lambda_{k}=-\left(\frac{2\pi k}{L}\right)^{2}$.

In the pseudo-spectral method the solution and the source term are
evaluated on a grid of points. Consider a uniform grid of $n+1$ points

\begin{align}
x_{j} & =j\frac{L}{n},\quad j=0\ldots n.\label{eq:periodic-grid}
\end{align}
Due to the periodic boundary conditions we can consider only $n$
points during the computations

\begin{align}
f_{n} & =f\left(x_{n}\right)=f(x_{0})=f_{0}.
\end{align}

The point values $f_{j}$ can be transformed using

\begin{align}
f_{j} & =\sum_{k=0}^{n-1}\hat{f}_{k}e^{2\pi ikj/n},\label{eq:IDFT}\\
\hat{f}_{k} & =\frac{1}{n}\sum_{j=0}^{n-1}f_{j}e^{-2\pi ikj/n}.\label{eq:DFT}
\end{align}
Eq. \eqref{eq:DFT} is the discrete Fourier transform (DFT)
which is used to compute the coefficients of the vector $f_{j}$ in
the basis of the discrete Fourier modes $e^{2\pi ikj/n}$. Eq.
\eqref{eq:IDFT} is the inverse discrete Fourier transform (IDFT) and
is used to compute the point values $f_{j}$.

The pseudo-spectral approximate solution can be computed using 
Eq. \eqref{eq:spectral-1d-solution} for the coefficients of
the discrete Fourier modes on the finite grid.

When the boundary conditions are not periodic, different set of basis
functions is used. For the Dirichlet boundary conditions

\begin{align}
f & =0\quad\mathrm{on}\:\partial\Omega
\end{align}
the Fourier series uses only the real sine functions

\begin{align}
f(x) & =\sum_{k=1}^{\infty}\widehat{f}_{k}\sin\left(\pi kx/L\right),\label{eq:Four-series-1}\\
\hat{f}_{k} & =\frac{2}{L}\int_{0}^{L}f(x)\sin\left(\pi kx/L\right)dx,
\end{align}
while for the Neumann boundary conditions

\begin{align}
\frac{\partial f}{\partial x} & =0\quad\mathrm{on}\:\partial\Omega
\end{align}
the Fourier series uses only the real cosine functions

\begin{align}
f(x) & =\frac{\widehat{f}_{0}}{2}+\sum_{k=1}^{\infty}\widehat{f}_{k}\cos\left(\pi kx/L\right),\label{eq:Four-series-1-1}\\
\hat{f}_{k} & =\frac{2}{L}\int_{0}^{L}f(x)\cos\left(\pi kx/L\right)dx,
\end{align}
and the Fourier coefficients $\hat{f}_{k}$ are real numbers.

Also, two different kinds of the uniform grid are possible. The first
type of the grid is the regular grid, which in the case of the Dirichlet
boundary conditions contains only the internal points

\begin{align}
x_{j} & =(j+1)\frac{L}{n+1},\quad j=0\ldots n-1,
\end{align}
whereas in the case of the Neumann boundary conditions it also contains
the boundary nodes

\begin{align}
x_{j} & =j\frac{L}{n-1},\quad j=0\ldots n-1.
\end{align}

The other type is the staggered grid

\begin{align}
x_{j} & =\left(j+\frac{1}{2}\right)\frac{L}{n},\quad j=0\ldots n-1.
\end{align}

In the psudo-spectral method on a grid different discrete sine transforms
(DST) and discrete cosine transforms (DCT) have to be used depending
on the grid and the boundary conditions. They are summarized in Table
\ref{tab:DST_DCT}.

\begin{table}
\caption{Real discrete Fourier transforms}
\label{tab:DST_DCT}

\centering{}%
\begin{tabular}{cc}
\toprule 
DST-I & $\hat{f}_{k}=2\sum_{j=0}^{n-1}f_{j}\sin\left(\pi(j+1)(k+1)/(n+1)\right)$\tabularnewline
\midrule 
DST-II & $\hat{f}_{k}=2\sum_{j=0}^{n-1}f_{j}\sin\left(\pi(j+1/2)(k+1)/n\right)$\tabularnewline
\midrule 
DST-III & $\hat{f}_{k}=(-1)^{k}f_{n-1}+2\sum_{j=0}^{n-2}f_{j}\sin\left(\pi(j+1)(k+1/2)/n\right)$\tabularnewline
\midrule 
DCT-I & $\hat{f}_{k}=f_{0}+(-1)^{k}f_{n-1}+2\sum_{j=1}^{n-2}f_{j}\cos\left(\pi jk/(n-1)\right)$\tabularnewline
\midrule 
DCT-II & $\hat{f}_{k}=2\sum_{j=0}^{n-1}f_{j}\cos\left(\pi(j+1/2)k/n\right)$\tabularnewline
\midrule 
DCT-III & $\hat{f}_{k}=f_{0}+2\sum_{j=1}^{n-1}f_{j}\cos\left(\pi j(k+1/2)/n\right)$\tabularnewline
\bottomrule
\end{tabular}
\end{table}

\begin{table}
\caption{Discrete Fourier transforms used in the Poisson solver}
\label{tab:tranforms_BC}

\centering{}%
\begin{tabular}{cccc}
\toprule 
boundary conditions & grid & forward & backward\tabularnewline
\midrule
\midrule 
periodic & regular & DFT & IDFT\tabularnewline
\midrule 
Dirichlet & regular & DST-I & $\frac{1}{2(n+1)}$DST-I\tabularnewline
\midrule 
Dirichlet & staggered & DST-II & $\frac{1}{2n}$DST-III\tabularnewline
\midrule 
Neumann & regular & DCT-I & $\frac{1}{2(n-1)}$DCT-I\tabularnewline
\midrule 
Neumann & staggered & DCT-II & $\frac{1}{2n}$DCT-III\tabularnewline
\bottomrule
\end{tabular}
\end{table}

Finally the eigenvalues are specified in the same notation as the
discrete transform in Table \ref{tab:spectral_eigenvalues}. The whole
algorithm of the Poisson solver can be summarized as
\begin{enumerate}
\item Compute the coefficients of the right-hand side $g$ in the basis
of eigenvectors of the discrete Laplacian using the forward discrete
transforms in Table \ref{tab:tranforms_BC}.
\item Compute the solution in this basis by dividing the coefficients by
the eigenvalues in Table \ref{tab:spectral_eigenvalues}.
\item Transform the solution back to the basis of the grid point values
using the backward transforms in Table \ref{tab:tranforms_BC}.
\end{enumerate}
The efficiency of this algorithm comes from the possibility of using
the Fast Fourier Transform (FFT) to perform the discrete transforms.

\begin{table}
\caption{The eigenvalues for the eigenvectors corresponding to the transforms
in Table \ref{tab:tranforms_BC} in the pseudo-spectral approximation }
\label{tab:spectral_eigenvalues}

\centering{}%
\begin{tabular}{ccl}
\toprule 
boundary conditions & grid & eigenvalues\tabularnewline
\midrule
\midrule 
\multirow{2}{*}{periodic} & \multirow{2}{*}{regular} & $\lambda_{k}=-\left(\frac{2\pi k}{L}\right)^{2},\quad k<[n/2]-1$\tabularnewline
 &  & $\lambda_{k}=-\left(\frac{2\pi(n-k)}{L}\right)^{2},\quad k\geqq[n/2]-1$\tabularnewline
\midrule 
Dirichlet & regular & \emph{$\lambda_{k}=-\left(\frac{\pi(k+1)}{L}\right)^{2}$}\tabularnewline
\midrule 
Dirichlet & staggered & $\lambda_{k}=-\left(\frac{\pi(k+1)}{L}\right)^{2}$\tabularnewline
\midrule 
Neumann & regular & $\lambda_{k}=-\left(\frac{\pi k}{L}\right)^{2}$\tabularnewline
\midrule 
Neumann & staggered & $\lambda_{k}=-\left(\frac{\pi k}{L}\right)^{2}$\tabularnewline
\bottomrule
\end{tabular}
\end{table}

\subsection{Finite difference method}

In certain application (e.g., computational fluid dynamics) it may be necessary to discretize 
the Poisson equation by a different scheme to ensure compatibility with other parts of a larger solver.
A common way to discretize Eq. \eqref{eq:poisson-continuous} are the second order central finite differences
\begin{align}
\frac{\partial^2\varphi}{\partial x^2}(x_{j}) & \simeq\frac{\varphi_{j+1}-2\varphi_{j}+\varphi_{j-1}}{\left(\triangle x\right)^{2}},\label{eq:FD2-Laplace}
\end{align}
where $\triangle x=x_{j}-x_{j-1}$. This scheme was used in most of the early studies on
fast Poisson solvers \citep{WilhelmsonEricksen, Swarztrauber:review, swarztrauber:1986symmetric}.

Discretization \eqref{eq:FD2-Laplace} leads to the
system of linear algebraic equations
\begin{align}
L\varphi & =g,\label{eq:FD2-matrix}
\end{align}
where $L$ is the matrix resulting from \eqref{eq:FD2-Laplace}
and from the boundary conditions. As reviewed by \citet{SchumannSweet:FFTs}
the eigenvectors of matrix $L$ are those used in the discrete transforms
in Table \ref{tab:tranforms_BC} with eigenvalues summarized in Table
\ref{tab:FD2-eigenvalues}.

\begin{table}
\caption{The eigenvalues for the eigenvectors corresponding to the transforms
in Table \ref{tab:tranforms_BC} in the second order central finite
difference approximation\label{tab:FD2-eigenvalues} according to
\citet{SchumannSweet:FFTs}}

\centering{}%
\begin{tabular}{ccl}
\toprule 
boundary conditions & grid & eigenvalues\tabularnewline
\midrule
\midrule 
\multirow{1}{*}{periodic} & \multirow{1}{*}{regular} & $\lambda_{k}=-\left(\frac{2\sin\left(\frac{k\pi}{n}\right)}{\Delta x}\right)^{2}$\tabularnewline
\midrule 
Dirichlet & regular & \emph{$\lambda_{k}=-\left(\frac{2\sin\left(\frac{\pi(k+1)}{2(n+1)}\right)}{\Delta x}\right)^{2}$}\tabularnewline
\midrule 
Dirichlet & staggered & $\lambda_{k}=-\left(\frac{2\sin\left(\frac{\pi(k+1)}{2n}\right)}{\Delta x}\right)^{2}$\tabularnewline
\midrule 
Neumann & regular & $\lambda_{k}=-\left(\frac{2\sin\left(\frac{\pi k}{2(n-1)}\right)}{\Delta x}\right)^{2}$\tabularnewline
\midrule 
Neumann & staggered & $\lambda_{k}=-\left(\frac{2\sin\left(\frac{\pi k}{2n}\right)}{\Delta x}\right)^{2}$\tabularnewline
\bottomrule
\end{tabular}
\end{table}

The algorithm for the solution of the algebraic system \eqref{eq:FD2-matrix}
is therefore equivalent to the algorithm for the pseudo-spectral approximation
in Section \ref{sub:Pseudo-spectral-method} with a different set
of eigenvalues.

\subsection{Multiple dimensions}

If the boundary conditions on all domain boundaries are equal, the
extension of the algorithm to multidimensional problems is straightforward.
The multidimensional discrete Fourier transforms are defined as separable
products of the one dimensional transforms. For example, for the periodic
boundary conditions in 3D the following transform can be used
\begin{align}
f_{j_{1}j_{2}j_{3}} & =\sum_{k_{1}=0}^{n_{1}-1}e^{2\pi ik_{1}j_{1}/n_{1}}\sum_{k_{2}=0}^{n_{2}-1}e^{2\pi ik_{2}j_{2}/n_{2}}\sum_{k_{3}=0}^{n_{3}-1}\hat{f}_{k_{1}k_{2}k_{3}}e^{2\pi ik_{3}j_{3}/n_{3}},\label{eq:IDFT-3D}\\
\hat{f}_{k_{1}k_{2}k_{3}} & =\frac{1}{n_{1}n_{2}n_{3}}\sum_{j_{1}=0}^{n_{1}-1}e^{-2\pi ik_{1}j_{1}/n_{1}}\sum_{j_{2}=0}^{n_{2}-1}e^{-2\pi ik_{2}j_{2}/n_{2}}\sum_{j_{3}=0}^{n_{3}-1}f_{j_{1}j_{2}j_{3}}e^{-2\pi ik_{3}j_{3}/n_{3}}.\label{eq:DFT-3D}
\end{align}

 Efficient multidimensional discrete transforms are
commonly included in the FFT libraries. Similar approach can be used
when the boundary conditions differ, as shown for the combination
of the periodic and the staggered Neumann boundary conditions by \citet{WilhelmsonEricksen}.
It is then necessary to use a sequence of transforms of smaller dimension.

The eigenvalues of the discrete 3D Laplace operator associated with
the Fourier mode $\hat{f}_{k_{1}k_{2}k_{3}}$ are computed from the
one dimensional eigenvalues as
\begin{align}
\lambda_{k_{1}k_{2}k_{3}} & =\lambda_{k_{1}}+\lambda_{k_{2}}+\lambda_{k_{3}}.
\end{align}

\section{Software implementation\label{sec:Software-implementation}}

The fast Poisson solver PoisFFT is a library written in Fortran 2003
with bindings to C and \CC. It uses the FFTW3\citep{FFTW}
library for the discrete Fourier transforms and the PFFT\citep{Pippig:PFFT}
library for the MPI parallelization of FFTW3 transforms. It is distributed
as a free software with the GNU GPLv3 license, which also covers FFTW3
and PFFT. The solver objects in PoisFFT are available in Fortran as
separate derived types in the IEEE single and double precision
and as a template class in \CC.

\begin{figure}

\begin{centering}
\includegraphics[height=4cm]{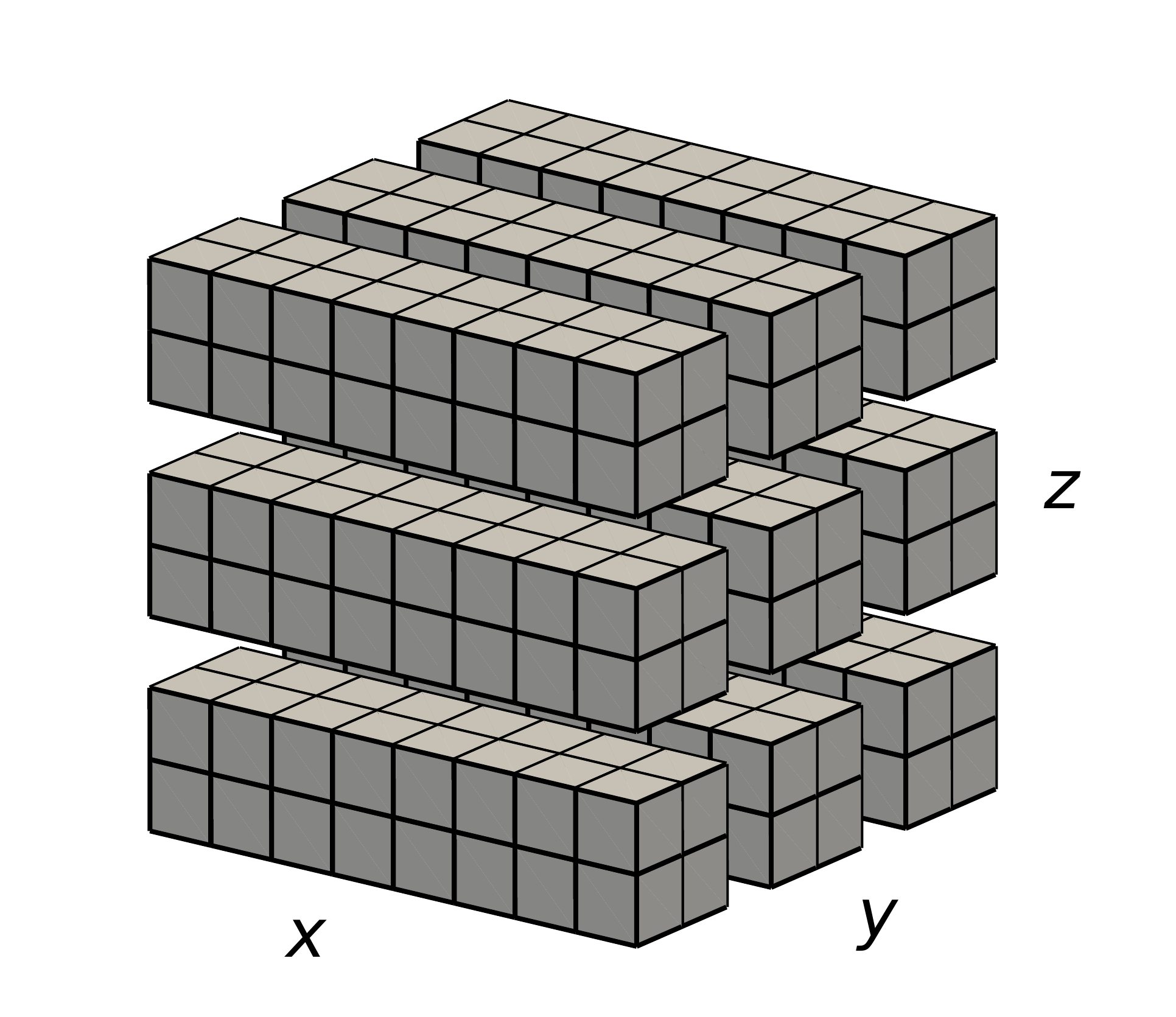}\caption{Decomposition of the domain in MPI into a 2D 3$\times$3 process grid.\label{fig:PoisFFT-2d-decomp}}

\par\end{centering}

\end{figure}

The library is parallelized using OpenMP for shared memory computers
and using MPI for distributed memory computers (e.g., clusters). Combination
of OpenMP and MPI at the same time is not implemented. It can also
employ the internal FFTW parallelization using POSIX threads. The
MPI version of the library uses the PFFT library to partition the
3D grid using a 2D decomposition along the last two dimensions in
Fortran (along the first two dimensions in C) as shown in Fig. \ref{fig:PoisFFT-2d-decomp}.

The main advantage of PoisFFT should be its simplicity and generality.
The initialization of the solver, the computation and the finalization of the solver are just 
three calls to the library. The arrays containing the right hand side or the solution
can be a subset of a larger array. This situation arises in computational fluid dynamics when ghost cells
are defined around the solution domain.

PoisFFT can solve the Poisson equation in 1, 2 and 3 dimensions with
any of the boundary conditions presented in Section \ref{sec:Solution-methods}
provided all domain boundaries use the same boundary condition. 

Because of practical applications in geophysical fluid dynamics, particularly
in the computational fluid dynamics model CLMM (Section \ref{sec:Application-to-incompressible}),
the combination of the periodic boundary condition in the $x$
direction and Dirichlet or Neumann boundary condition in the $z$ direction
is also implemented. In this case boundary condition in the $y$ direction can be periodic or
 identical to the condition used in the $z$ direction. 

Because the parallel real 1D transform is
not available in FFTW3 or PFFT we parallelized this step using
a global array transpose via the subroutine\texttt{ MPI\_Alltoallv()}
and local array transposes. This is necessary for the combined boundary conditions.

\section{Performance evaluation\label{sec:Performance-evaluation}}

Due to the differences in the discrete Fourier transforms for different
boundary conditions, the wall-clock times and the CPU times necessary
for the solution of the Poisson equation depend on the boundary conditions.
We chose the 3D solvers with the periodic boundary conditions and
the Neumann boundary conditions on the staggered grid to test the
solver with the 3D complex transform and with the 3D real transform.

All tests were performed at the cluster \texttt{zapat.cerit-sc.cz}
at the Czech supercomputing center Cerit-SC. The computational nodes
contained 2 sockets with 8-core Intel E5-2670 2.6 GHz processors and
128 GB RAM and were connected using an InfiniBand network. PoisFFT
was compiled with the Intel Fortran compiler 14.1 and linked with
OpenMPI 1.8.

\begin{figure}
\begin{centering}
\begin{tabular}{cc}
\includegraphics[width=0.5\textwidth]{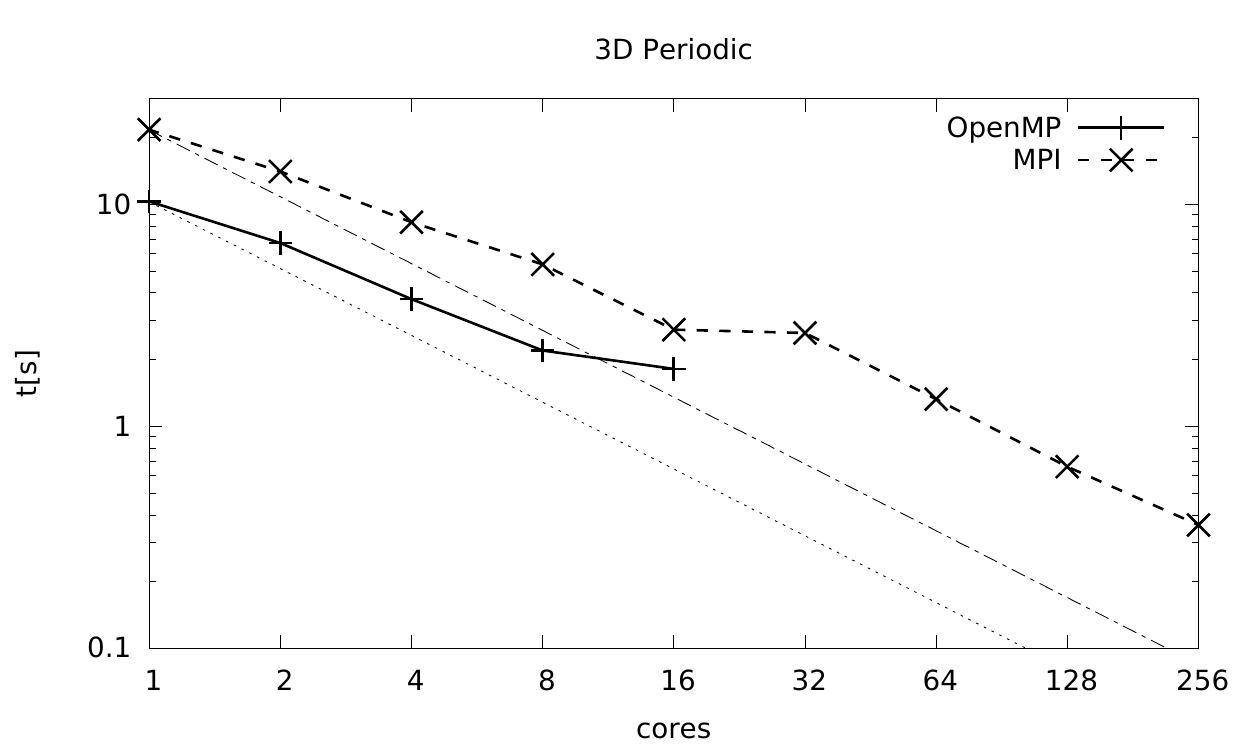} & \includegraphics[width=0.5\textwidth]{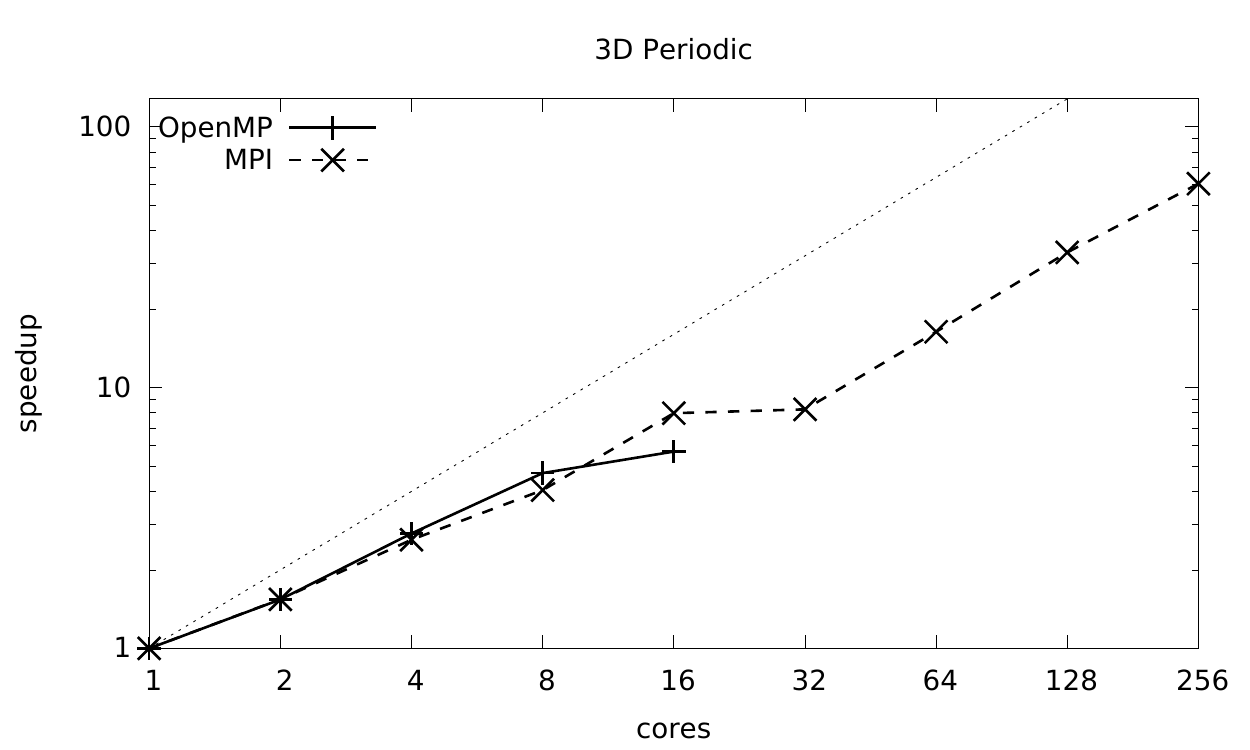}\tabularnewline
\includegraphics[width=0.5\textwidth]{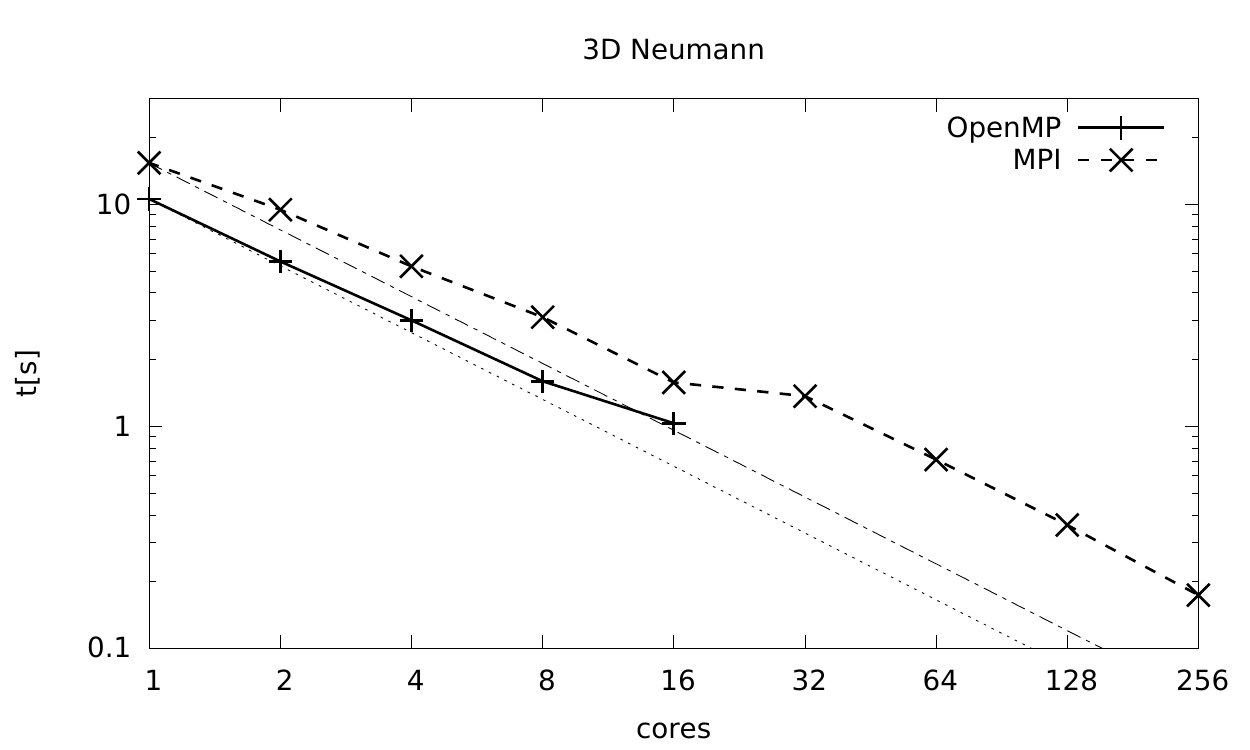} & \includegraphics[width=0.5\textwidth]{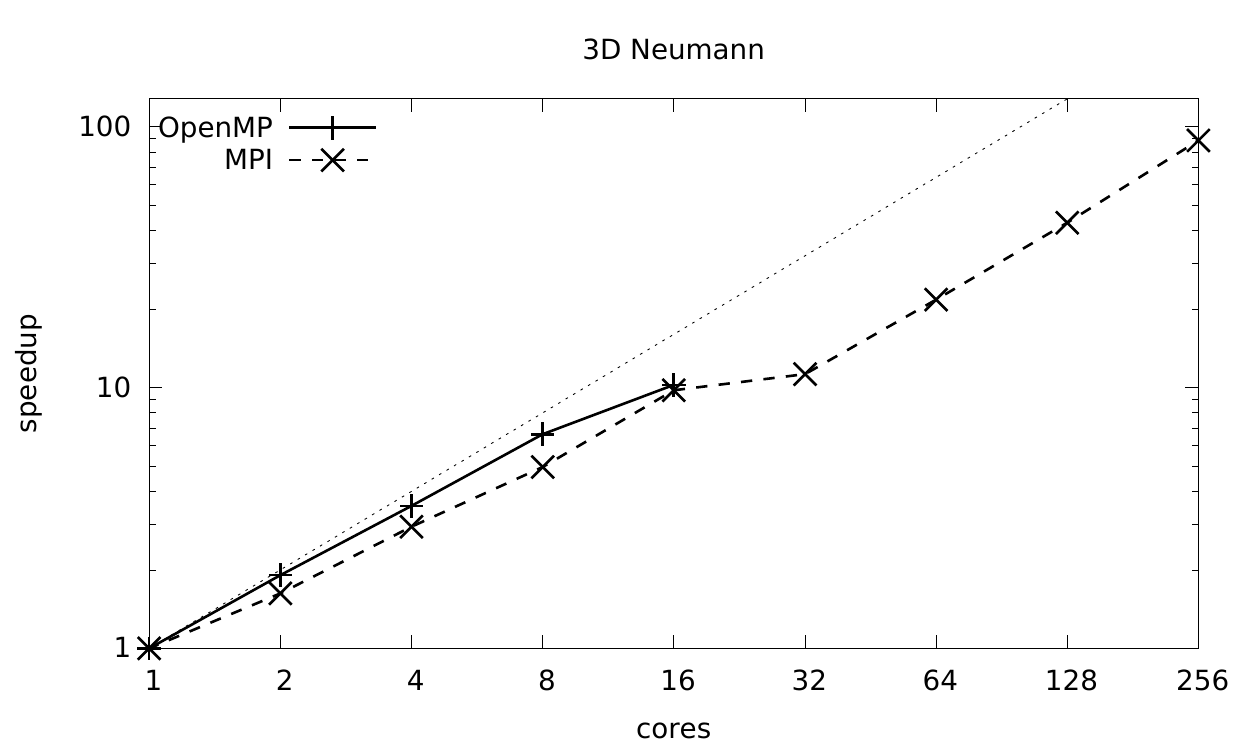}\tabularnewline
\end{tabular}
\par\end{centering}

\caption{Strong scaling of PosFFT with $512^{3}$ points.\label{fig:PoisFFT-strong-scaling}}

\end{figure}

The first test examined the strong scaling on a fixed domain with
512 points in every dimension. The same problem was solved with an
increasing number of CPU cores. From the measured times in Fig. \ref{fig:PoisFFT-strong-scaling}
it is clear that the OpenMP version, which directly calls FFTW, is
considerably faster than the MPI version, which calls the PFFT wrapper
over FFTW and which has to perform significant amount of message passing
between processes. However, OpenMP is limited to single node on the
cluster used for the computations. One can also notice a decrease
of efficiency when going from 16 to 32 cores. One node has 16 cores
and the expensive network memory transfers slow the computation down when
more nodes are needed.

\begin{figure}

\begin{centering}
\begin{tabular}{cc}
\includegraphics[width=0.5\textwidth]{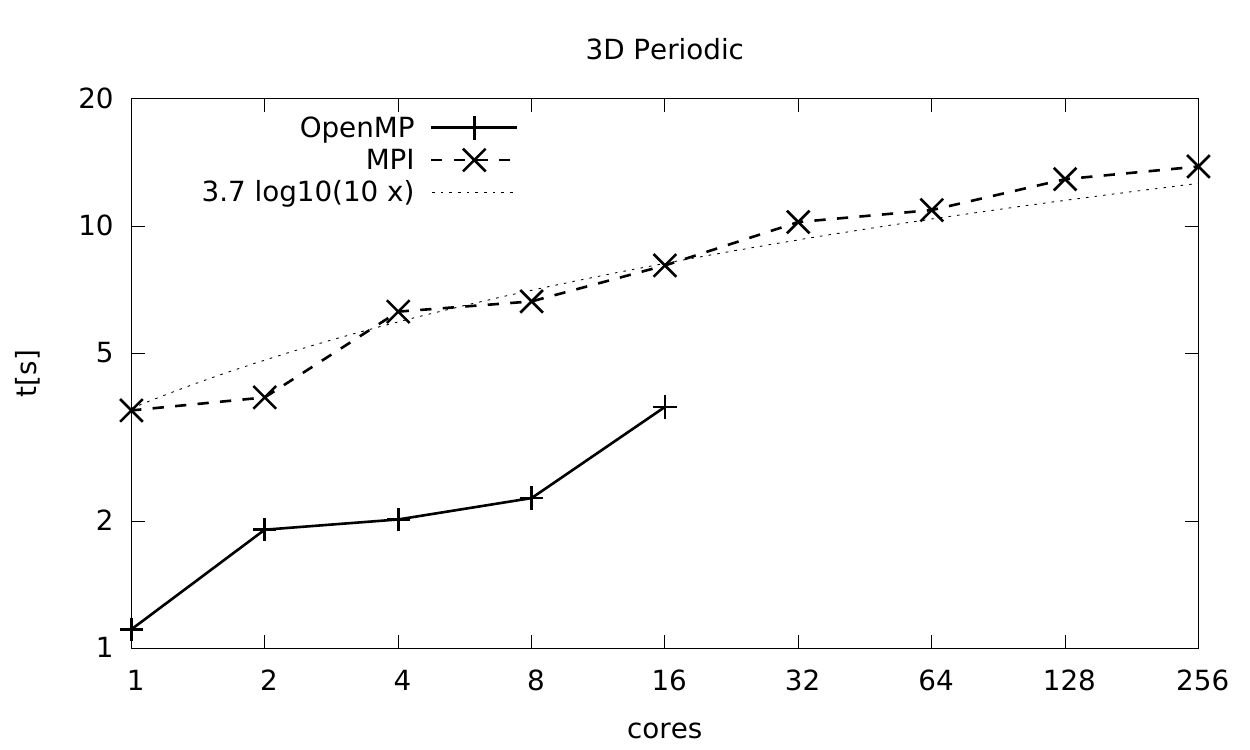} & \includegraphics[width=0.5\textwidth]{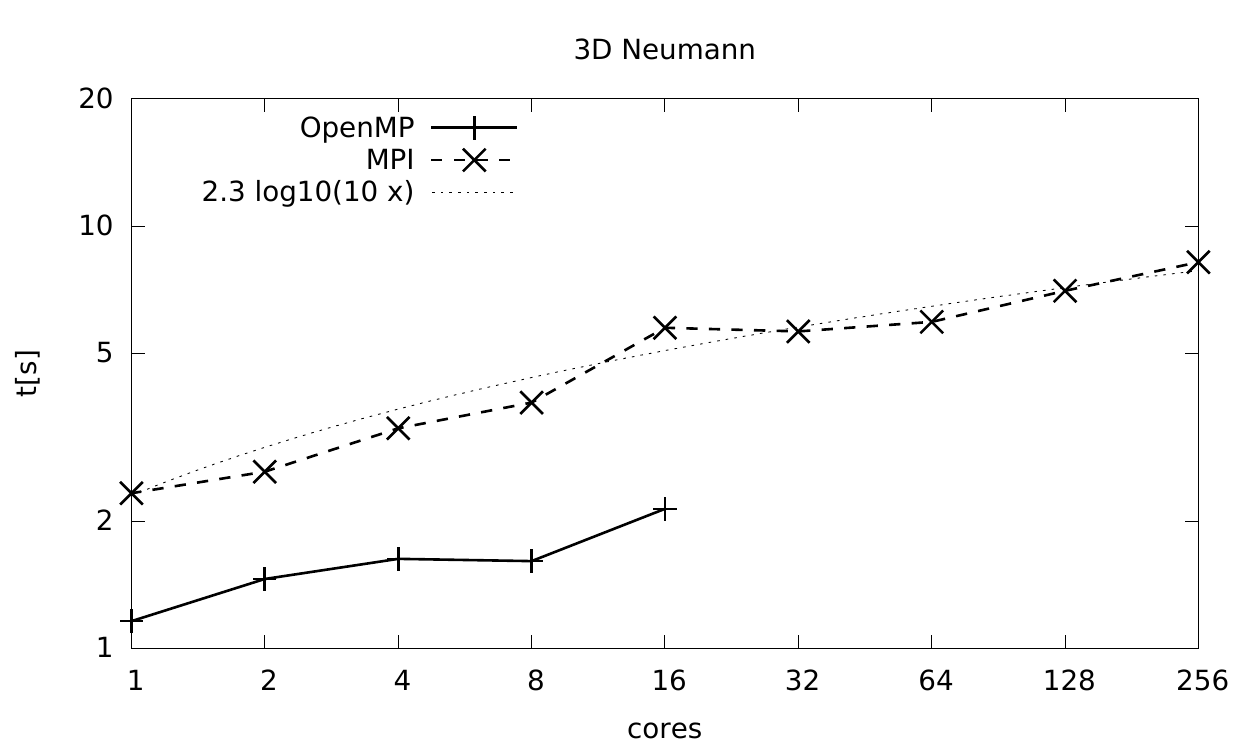}\tabularnewline
\end{tabular}\caption{Weak scaling of PoisFFT with $256^{3}$ points per core.\label{fig:PoisFFT-scaling-weak}}

\par\end{centering}

\end{figure}

For the weak scaling test we used a domain discretization in which
the number of points increased with the number of CPU cores. There
were 256 points in each dimension per core. The complete problem size
was therefore proportional to the number of cores. The results of
the time measurements are in Fig. \ref{fig:PoisFFT-scaling-weak}.
They are consistent with the scaling $\mathcal{O}\left(N\ln\left(N\right)\right)$
which holds for the FFT transforms\citep{FFTW}.

\section{Application of PoisFFT to simulation of incompressible flow\label{sec:Application-to-incompressible}}

\subsection{Model CLMM}

Originally, PoisFFT was developed to support the atmospheric computational
fluid dynamics code CLMM (Charles University Large-Eddy Microscale
Model)\citep{fuka:FVCA}. This model simulates atmospheric flows using
the incompressible Navier-Stokes equations with an eddy viscosity
subgrid model.

The incompressible Navier-Stokes equations are

\begin{align}
\frac{\partial\boldsymbol{u}}{\partial t}+\nabla\cdot(\boldsymbol{u}\otimes\boldsymbol{u}) & =-\nabla p+\nu_{\mathrm{sgs}}\nabla^{2}\boldsymbol{u},\label{eq:pois-incomp-mom}\\
\nabla\cdot\boldsymbol{u} & =0,\label{eq:pois-incomp-cont}
\end{align}
where $\boldsymbol{u}$ is the velocity vector, $p$ is the pressure
and $\nu_{\mathrm{sgs}}$ is the eddy viscosity modeled according
to \citet{nicoud:sigma}. Eqs. (\ref{eq:pois-incomp-mom})--(\ref{eq:pois-incomp-cont})
do not contain an evolution equation for pressure. The model CLMM
uses the projection method\citep{brown:proj} and the 3 stage Runge-Kutta
method\citep{spalart:rk3}

\begin{align}
\boldsymbol{u}^{(1)} & =\boldsymbol{u}^{n},\: p^{(1)}=p^{n-1/2},\label{eq:pois-RK3expl-begin}\\
\frac{\boldsymbol{u}^{*}-\boldsymbol{u}^{(k)}}{\Delta t} & =-\alpha_{k}\nabla p^{(k)}+\gamma_{k}\left(-\nabla\cdot(\boldsymbol{u}^{(k)}\boldsymbol{u}^{(k)})+\nu\nabla^{2}\boldsymbol{u}^{(k)}\right)+\nonumber \\
 & \qquad+\zeta_{k}\left(-\nabla\cdot(\boldsymbol{u}^{(k-1)}\boldsymbol{u}^{(k-1)})+\nu\nabla^{2}\boldsymbol{u}^{(k-1)}\right)\label{eq:pois-RK3expl-mom}\\
\nabla^{2}\phi & =\frac{\nabla\cdot\boldsymbol{u}^{*}}{\alpha_{k}\Delta t},\label{eq:projection-Poisson}\\
\frac{\boldsymbol{u}^{(k+1)}-\boldsymbol{u}^{*}}{\alpha_{k}\Delta t} & =-\nabla\phi,\\
p^{(k+1)} & =p^{(k)}+\phi,\label{eq:pois-RK3expl-prupdate}\\
\mathrm{for}\; & k=1..3,\quad\mathrm{and}\nonumber \\
\boldsymbol{u}^{n+1} & =\boldsymbol{u}^{(4)},\: p^{n+1/2}=p^{(4)},\label{eq:pois-RK3expl-end}
\end{align}
with the set of coefficients

\begin{align}
\alpha_{1} & =\frac{8}{15},\:\alpha_{2}=\frac{2}{15},\:\alpha_{3}=\frac{1}{3},\nonumber \\
\gamma_{1} & =\frac{8}{15},\:\gamma_{2}=\frac{5}{12},\:\gamma_{3}=\frac{3}{4},\nonumber \\
\zeta_{1} & =0,\:\zeta_{2}=-\frac{17}{60},\:\zeta_{3}=-\frac{5}{12}.\label{eq:RK3-coefs}
\end{align}

The spatial derivatives in (\ref{eq:pois-RK3expl-begin})--(\ref{eq:pois-RK3expl-end})
are discretized using the second order central finite differences
on a staggered grid\citep{KimKimChoi} with the exception of the advective
terms $-\nabla\cdot(\boldsymbol{u}\otimes\boldsymbol{u})$ for which fourth order
central differences are used.

In every Runge-Kutta stage first a predictor velocity field $\boldsymbol{u}^{*}$
is computed and afterward it is corrected to be divergence free. The
discrete Poisson problem \eqref{eq:projection-Poisson} must be solved
in the corrector step, with the staggered Neumann boundary conditions
on all outside boundaries except the boundaries where the periodic
boundary condition for the velocity are imposed and where the periodic
condition applies to the discrete Poisson equation as well. It is the task
of the PoisFFT library to solve Eq. \eqref{eq:projection-Poisson}.

It is possible to use this model on the uniform grid to simulate certain types of complex
geometries thanks to the immersed boundary method (IBM)\citep{IBM:review},
which introduces additional sources of momentum and mass inside or
near the solid obstacles. In CLMM the IBM interpolations of \citet{Peller_etal:IBM}
and the mass sources according to \citet{KimKimChoi} are implemented.

\subsection{Performance test}

To evaluate the effect of CLMM we chose the test case defined by the
COST Action ES1006 ``Evaluation, improvement and guidance for the
use of local-scale emergency prediction and response tools for airborne
hazards in built environments'' with a virtual city called Michelstadt
which serves for validation of models for flow and dispersion in urban
environments\citep{CEDVAL-LES}. The description of the geometry and
the flow validation data are available at the CEDVAL-LES website%
\footnote{\texttt{http://www.mi.uni-hamburg.de/CEDVAL-LES-V.6332.0.html, accessed
on September 28th, 2014.}%
}. In this paper only the velocity field is simulated for a direct
comparison of the relative performance of the Poisson solver PoisFFT
and other parts of the model.

The computational domain of the Michelstadt city had dimensions 1326
m $\times$ 836 m $\times$ 399 m and was discretized using a constant
resolution of 3m in all directions leading to the total number of
16.3 million cells. The performance was measured for the first 100
time-steps for different parallel setups -- OpenMP and MPI -- with
different numbers of CPU cores. The fastest configuration was used
to compute the flow for 10000 s of the real scale time and the average
flow and turbulence was compared to the experimental values. The performance
was tested on the cluster \texttt{luna.fzu.cz} in the Czech national
grid MetaCentrum. The computational nodes have 2 sockets with the
8-core Intel Xeon E5-2650 2.6 GHz CPU, 96 GB of RAM and are connected
using an InfiniBand network.

\begin{figure}
\begin{centering}
\includegraphics[width=0.5\textwidth]{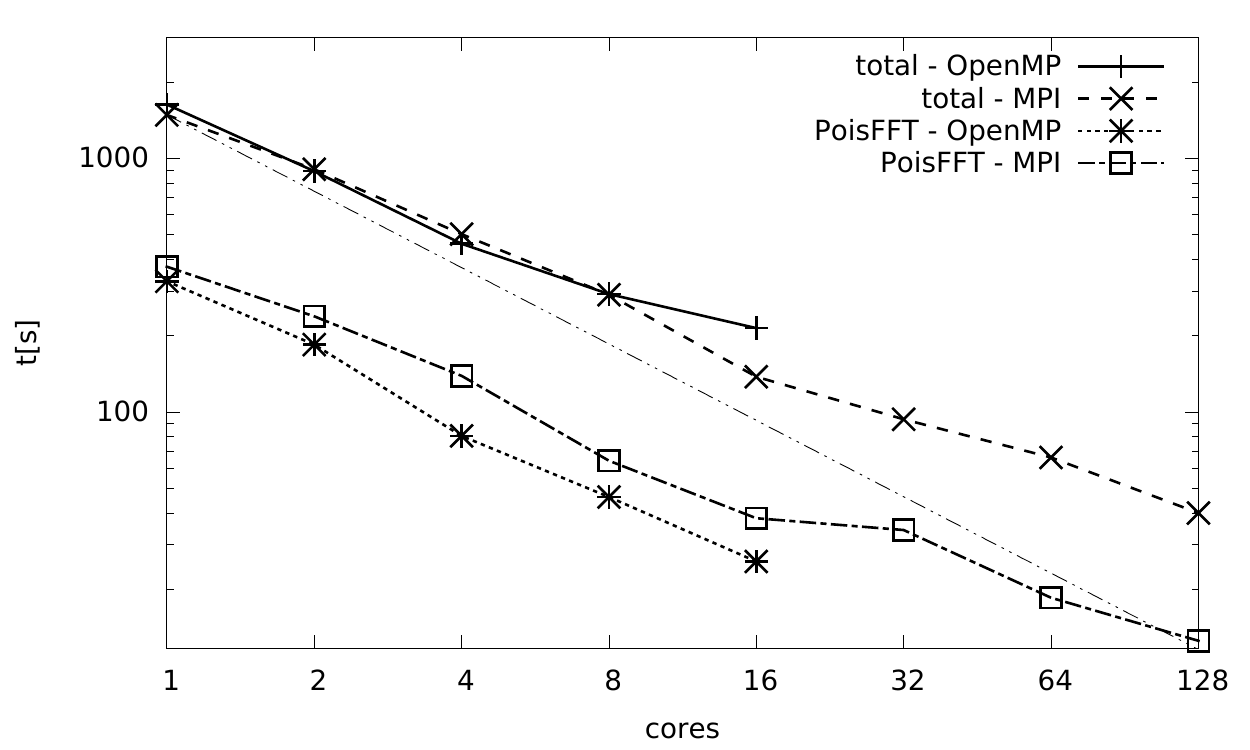}\includegraphics[width=0.5\textwidth]{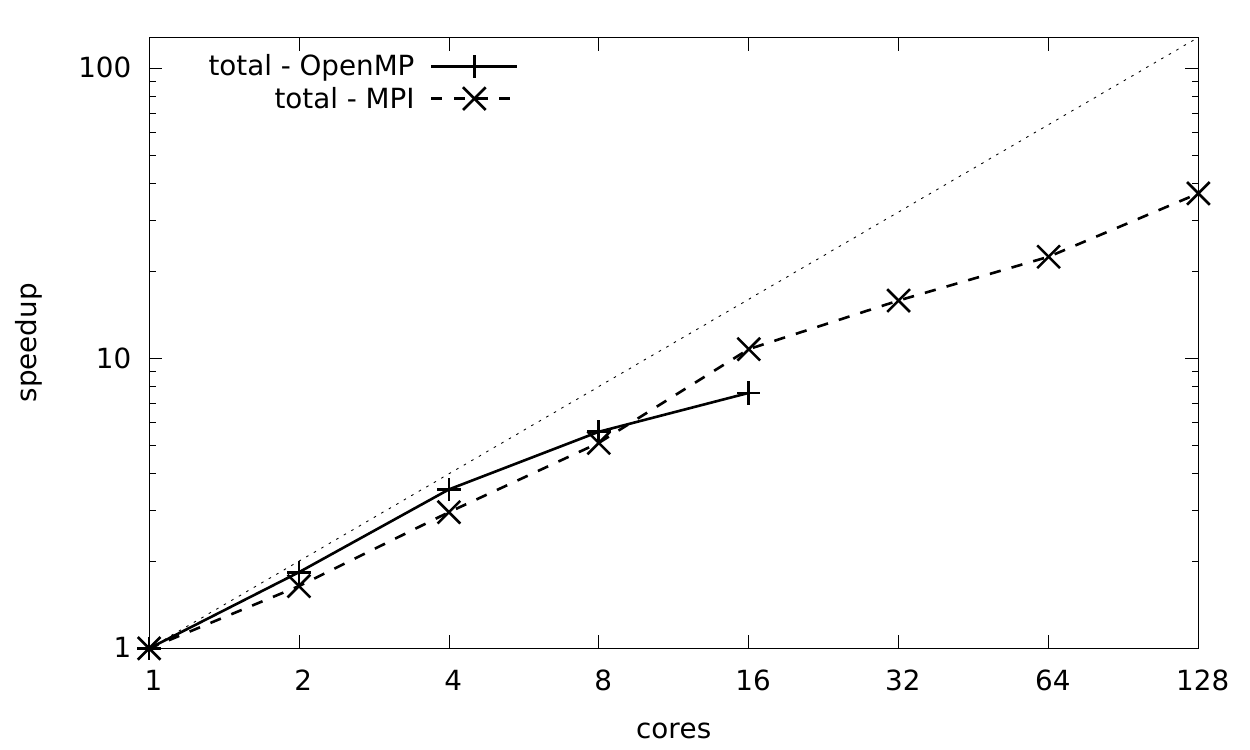}
\par\end{centering}

\caption{Strong scaling of CLMM and PoisFFT.\label{fig:CLMM-strong-scaling}}

\end{figure}

The scaling of the wall-clock time for the solution of first 100 time-steps
is in Fig. \ref{fig:CLMM-strong-scaling}. Although the MPI version
of PoisFFT is slower than the OpenMP version for the same number of
CPU cores, the whole model CLMM is as fast with MPI as with OpenMP
on the used nodes up to 8 threads. The OpenMP version with 16 threads
is slower than the MPI version with 16 processes. CLMM appears
to benefit from the improved data locality in the
distributed MPI version. The relative contribution of PoisFFT to the
total CLMM solution time is between 22 \% -- 36 \% when using MPI.
When using OpenMP the PoisFFT portion starts at 20 \% for one and
two threads and decreases for higher numbers of threads as the parallel
efficiency of CLMM in OpenMP decreases. From these ratios it is clear
that the efficient Poisson solver is a necessity for an efficient
computation and that PoisFFT scales well enough so that it doesn't
become a bottleneck at the moderate numbers of cores used for the
test.

\begin{figure}
\begin{centering}
\begin{tabular}{cc}
\includegraphics[width=0.48\textwidth]{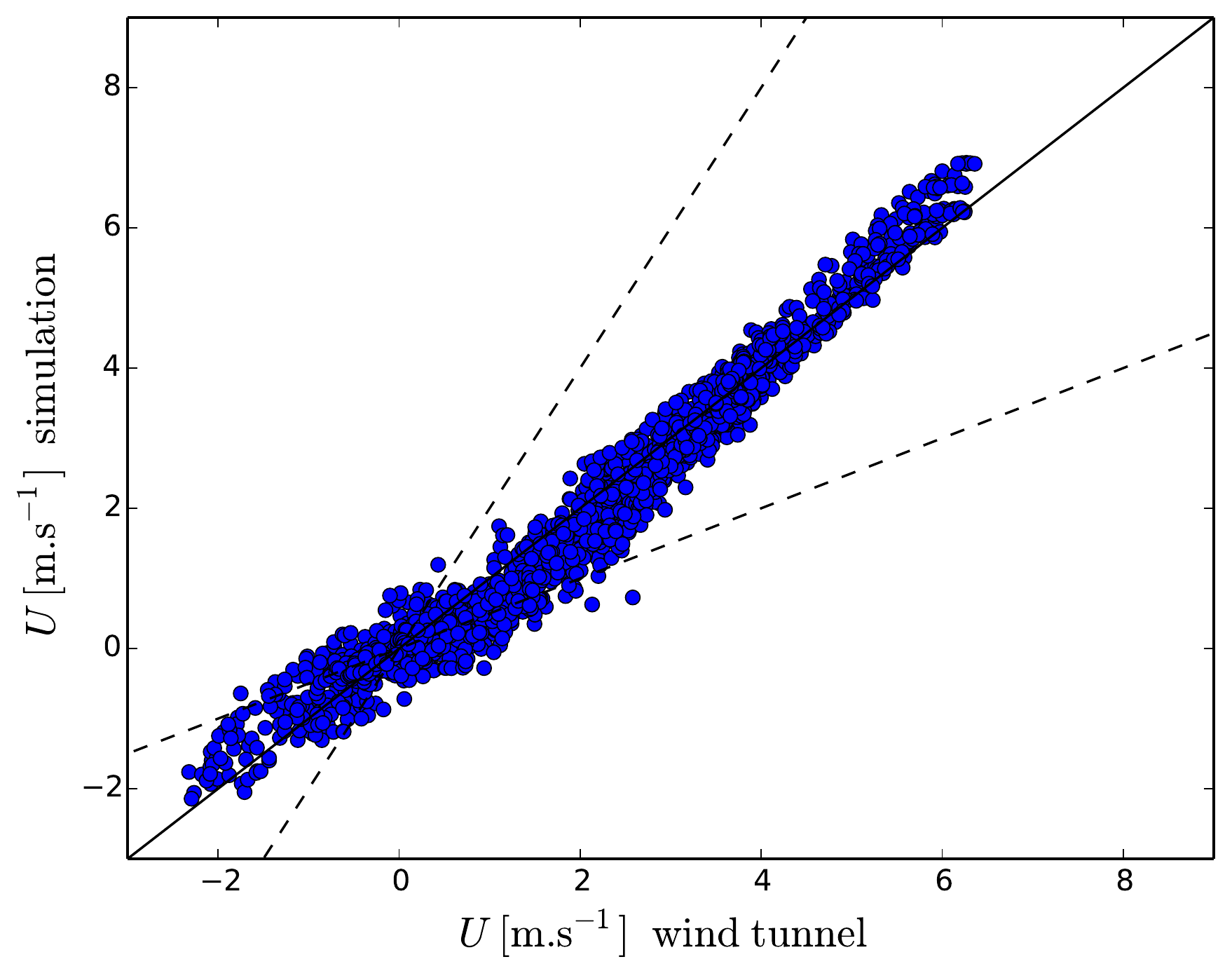} & \includegraphics[width=0.48\textwidth]{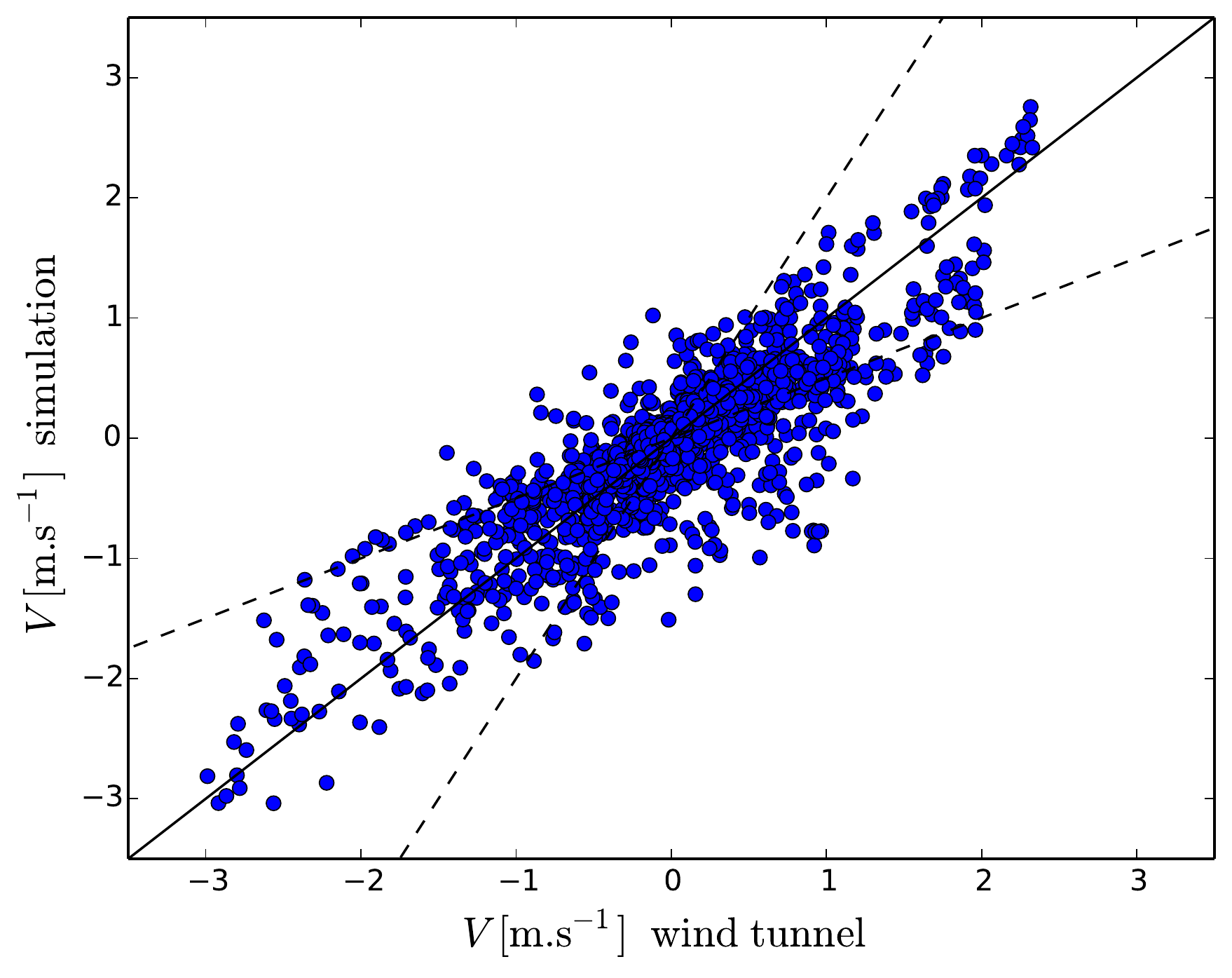}\tabularnewline
\end{tabular}
\par\end{centering}

\caption{Scatter plots for the comparison of simulated and measured mean stream-wise
($U$) and lateral ($V$) velocity components.\label{fig:PoisFFT-Michelstadt-scatters}}
\end{figure}

\begin{figure}
\includegraphics[width=1\textwidth]{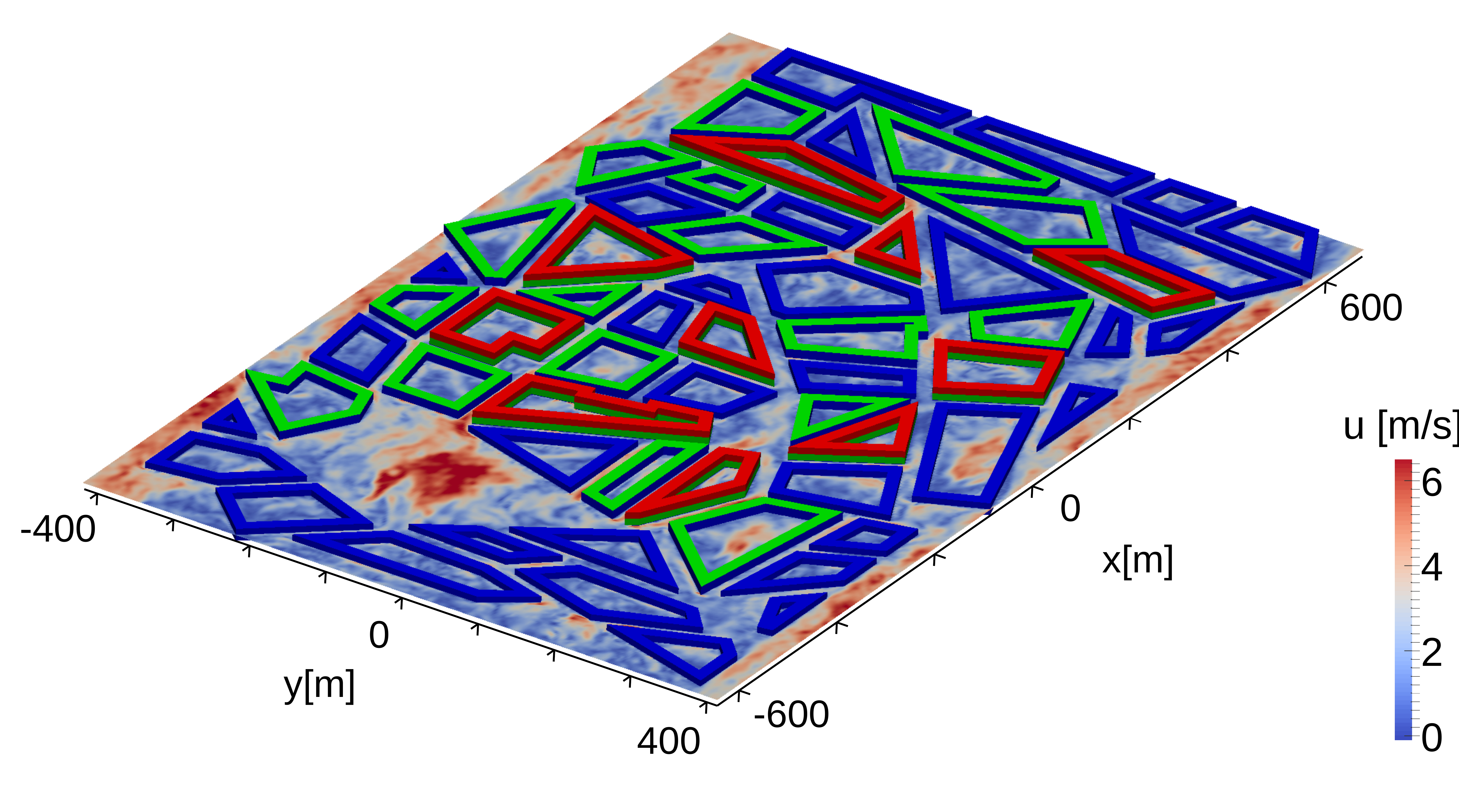}

\caption{The Michelstadt geometry with the instantaneous field of the wind
velocity magnitude at height 7.5 m above ground. The blue, green and
red buildings are 15, 18 and 24 m high, respectively.\label{fig:PoisFFT-Michelstadt-field}}

\end{figure}

The correctness of the computation was evaluated by comparison of
the time-averaged simulated wind field with the measurements. Fig.
\ref{fig:PoisFFT-Michelstadt-scatters} shows the comparison of the
values of the velocity components in the reference points of CEDVAL-LES.
An illustration of the simulated wind field at height 7.5 m above
the ground is in Fig. \ref{fig:PoisFFT-Michelstadt-field}.

\section{Conclusions}

The fast Poisson solver PoisFFT is able to compute a discrete approximate
solution to the Poisson equation in the pseudo-spectral or second
order finite difference approximation. It is parallel using OpenMP
and MPI and the parallel scaling was demonstrated to be usable for
moderate numbers of CPU cores. The scaling for the number of cores
in the order of thousands remain untested due to the lack of suitable
hardware. The experiments of \citet{Pippig:PFFT} with the underlying
PFFT library suggest that the solver could scale even for thousands or 
tens of thousands of cores.

It was demonstrated that PoisFFT can be employed for the solution
of the discrete Poisson equation in the projection step in an incompressible
flow solver. It does not become a bottleneck at 128 cores and it is
expected to scale well even further.

The library can be used from Fortran, C and \CC. Interfaces to other
computer languages (e.g., Python) are planned.

\section*{Acknowledgment}

Computational resources were provided by the MetaCentrum under the
program LM2010005 and the CERIT-SC under the program \foreignlanguage{british}{Centre}
CERIT Scientific Cloud, part of the Operational Program Research and
Development for Innovations, Reg. no. CZ.1.05/3.2.00/08.0144. The
Michelstadt flow simulation was done in the framework of the COST
Action ES1006. The work was supported by the University Research Centre
UNCE 204020/2012.

\bibliographystyle{elsarticle-num-names}
\bibliography{moje}

\end{document}